\documentclass[12pt,a4paper]{amsart}

\usepackage{amsmath, amsfonts, xifthen, latexsym, amssymb, amsthm, amscd}

\usepackage[utf8]{inputenc}
\usepackage[shortlabels]{enumitem}
\usepackage{graphicx}
\usepackage{url}

\usepackage{hyperref}
\hypersetup{colorlinks=true,citecolor=blue,filecolor=blue,linkcolor=blue,urlcolor=blue}
\usepackage[margin=1.4in]{geometry}

\newtheorem{theorem}{Theorem}
\newtheorem{lemma}{Lemma}[section]

\newtheorem{proposition}{Proposition}[section]

\theoremstyle{definition}
\newtheorem{definition}{Definition}[section]

\newtheorem{remark}{Remark}[section]

\renewcommand{\phi}{{\varphi}}

\newcommand\R{\mathbb R}

\newcommand\Z{\mathbb Z}

\newcommand\N{\mathbb N}

\newcommand{\actson}{\curvearrowright}

\newcommand{\cal}[1]{{\mathcal #1}}

\usepackage{etoolbox}
\newtoggle{no_cases}\toggletrue{no_cases}
\newcommand{\case}[2][]{\iftoggle{no_cases}{\left\{\begin{array}{ll}#2 & #1}{\\#2 & #1}\togglefalse{no_cases}}
\newcommand{\esac}{\end{array}\right.\toggletrue{no_cases}}

\newcommand{\Sub}{\operatorname{Sub}}

\newcommand{\K} {{\bf K}^*}
\begin{document}
\title[Amenable purely infinite actions]{Amenable purely infinite actions
on the non-compact Cantor set }
\author{Gábor Elek}
\address{Department of Mathematics And Statistics, Fylde College, Lancaster University, Lancaster, LA1 4YF, United Kingdom}

\email{g.elek@lancaster.ac.uk}  

\thanks{The author was partially supported
by the ERC Consolidator Grant "Asymptotic invariants of discrete groups,
sparse graphs and locally symmetric spaces" No. 648017. }

\begin{abstract} We prove that any countable non-amenable
group $\Gamma$ admits
a free minimal amenable purely infinite action on the 
non-compact Cantor set. This answers a question of Kellerhals, Monod and
R{\o}rdam \cite{KMR}.

\end{abstract}\maketitle
\tableofcontents
\newpage

\section{Introduction}
Free, minimal (topological) actions on the non-compact Cantor set $\K$ were constructed
by Matui and R{\o}rdam \cite{MR} and Danilenko \cite{Dan2} (see also \cite{Dan1}). Furthermore, the action constructed
in \cite{Dan2} is amenable (even Borel-hyperfinite). In \cite{KMR} Kellerhals, Monod and R{\o}rdam 
studied free minimal amenable actions on $\K$ that are purely infinite as well. They showed that
for such actions the associated reduced $C^*$-algebra is always a stable Kirchberg algebra of
the UCT class. Let us recall the notion of a  purely infinite action on $\K$ from \cite{KMR}
[Definition 4.4]. Let $\alpha:\Gamma\actson \K$ be an action of a 
countable group
$\Gamma$ on the non-compact Cantor set. Recall that $\K$ is the unique (up to homeomorphisms) locally compact,
 non-compact, totally disconnected, metrizable
Hausdorff space that contains no isolated points. We say that
a compact-open set $K$ is paradoxical with respect to the action if there exist pairwise disjoint
compact-open sets $K_1, K_2, \dots, K_{n+m}$  and elements $t_1,t_2,\dots,t_{n+m}\in\Gamma$
such that $K_j\subset K$ for all $j$ and
$$K=\cup^n_{i=1} \alpha(t_i)(K_i)=\cup^{n+m}_{j=n+1} \alpha(t_j)(K_j)\,.$$
\noindent
The action $\alpha$ is called {\bf purely infinite} if all the compact-open subsets of $\K$ are
paradoxical with respect to the action. In \cite{KMR} Kellerhals, Monod and R{\o}rdam  proved that 
if a countable group $\Gamma$ contains an exact non-supramenable subgroup, then $\Gamma$ admits
a free minimal amenable purely infinite action on $\K$. 
They asked whether any non-supramenable group admits such an action.
The  goal of this paper is to give a positive answer for this question by proving the following theorem.
\begin{theorem} \label{fotetel}
Every countable non-amenable group $\Gamma$ admits a free minimal amenable purely infinite action on $\K$.
\end{theorem}
\noindent
Note that Theorem \ref{fotetel} combined with the above mentioned result in 
\cite{KMR} implies that a countable group admits a free minimal
amenable purely infinite action on $\K$ if and only if it is non-supramenable.
We will prove the theorem for finitely generated groups and use the fact (Lemma 7.1 in \cite{KMR})  
that if a group $\Gamma$ contains a subgroup $H$ which
admits such free minimal amenable purely infinite action, then the group $\Gamma$ 
itself admits an action with the same properties as well. By the same reason, we do not need to consider the case of non-supramenable amenable groups.
The strategy of the proof goes as follows.
In Section \ref{noncom} we introduce the notion of non-compact Bernoulli shifts. Some special
elements of these spaces are called proper landscapes. In Section \ref{propland} we show that the
orbit closure of a proper landscape always contains an invariant
subset $Y$ that is homeomorphic to $\K$ and the Bernoulli action of the group on $Y$ is both free
and minimal. It is not hard to construct landscapes such that the resulting free minimal action is
Borel hyperfinite, but these actions cannot be extended to purely infinite actions. 
So, we will construct 
free minimal actions $\beta$  on $\K$ such that $\K$ can be exhausted by compacta that admit free group
actions from the topological full group of $\beta$.
Using the fact that the free group has Yu's
Property A we construct an amenable extension of $\beta$. Then, with the help 
of the partial actions of the free group and the well-known
paradoxical property of non-amenable graphs we inductively construct a sequence of extensions for which
more and more compact-open sets are becoming paradoxical, in such a way that freeness, minimality
and amenability are preserved. In the resulting free minimal amenable limit action all the compact-open sets will be paradoxical 
and this will finish the proof of our theorem. 

\section{Non-compact Bernoulli subshifts} \label{noncom}
 For the rest of
the paper let $\Gamma$ be a finitely generated group with a symmetric generating set $\Sigma=\{\sigma_i\}^n_{i=1}$. 
Let $A=\{0,1\}^\N\times \{\N\cup \{\infty\}\}$. We equip $\{0,1\}^\N$ with the standard product topology and we regard the space 
$\{\N\cup \{\infty\}\}$
as the compactification of the natural numbers. Hence, $A$ is a totally disconnected space homeomorphic to the Cantor set.
We consider the Bernoulli space $A^\Gamma$. Clearly, $A^\Gamma$ is homeomorphic to the Cantor set as well and $\Gamma$ acts (on the left) 
on $A^\Gamma$ continuously
by translations, that is, $L_\gamma(x)(\delta)=x(\delta\gamma)$.  If $x\in A^\Gamma$ and $x(\gamma)=(a,b)$, then we refer to $a$ as the 
Cantor coordinate $C(x(\gamma))$  of $x(\gamma)$
 and to $b$ as the height coordinate $H(x(\gamma))$ of $x(\gamma)$.
A minimal non-compact Bernoulli subshift is 
\begin{itemize}
\item a $\Gamma$-invariant subset $Y\subset A^\Gamma$ such that $Y$ is homeomorphic to $K^*$;
\item and for every element $y\in Y$
the orbit of $y$ is dense in $Y$.
\end{itemize}
\noindent
An element $y\in A^\Gamma$ is of  {\bf totally finite height} if for all $\gamma\in\Gamma$ the height coordinate of $y(\gamma)$ is a finite number. We say that
$y\in A^\Gamma$ is of {\bf totally infinite height} if for all $\gamma\in\Gamma$ the height coordinate of $y(\gamma)$ is $\infty$.
We call $y\in A^\Gamma$ of totally finite height a {\bf regular element} if the orbit closure of $y$ consists only of elements
of totally finite height and of totally infinite height. 
Our first goal is to give a sufficient condition for a regular element $y\in A^\Gamma$ having such orbit closure $\overline{O}(y)$ that the
totally finite height part of $\overline{O}(y)$  is a free, minimal, non-compact Bernoulli subshift. 

\noindent
Let us consider the (right) Cayley graph $G=\mbox{Cay}(\Gamma, \Sigma)$ of our group $\Gamma$ equipped with the shortest distance metric $d_G$.
A labeling $\lambda:\Gamma\to \{0,1\}^\N$ is called a {\bf proper Cantor labeling} if the following condition holds. For every $r>0$ there 
exists $S_r>0$
such that if $0<d_G(\gamma,\delta)\leq r$ then $(\lambda(\gamma))_{S_r}\neq (\lambda(\delta))_{S_r}$, where $(x)_s\in \{0,1\}^s$ denotes
the first $s$ coordinates of the element $x\in\{0,1\}^\N$. 
\begin{proposition}
There exist proper labelings on $\Gamma$.
\end{proposition}
\proof
Let $d$ be the degree of the vertices of $G$. Then, for all $k\geq 1$ we have a function
$$\lambda_k:\Gamma\to \{1,2,\dots, d^k+1\}$$
\noindent
such that if $0<d_G(\gamma,\delta)\leq k$ then
$$\lambda_k(\gamma)\neq \lambda_k(\delta)\,.$$
\noindent
Let
$\zeta_k:\Gamma\to \{0,1\}^{d^k+1}$ be
defined by
$$\zeta_k(\gamma)=(0,0,\dots,1,\dots,0,0)\,,$$
\noindent
where the only $1$ is at the $\lambda_k(\gamma)$-th position.
Now we can define the proper labeling by
$$\lambda(\gamma)=( \zeta_1(\gamma) \zeta_2(\gamma)\dots )\,.\quad \qed$$
\begin{lemma}
Let $y\in A^\Gamma$ be an element such that the Cantor coordinates of $y$ amount to a proper Cantor labeling of $\Gamma$. 
Then the action of
$\Gamma$ on the orbit closure of $y$ is free.
\end{lemma}
\proof
 Let $x\in \overline{O}(y)$. Then the Cantor
coordinates of $x$ also amount to a proper Cantor labeling of $\Gamma$. Consequently, the Cantor coordinates of $x$ are all different. 
Hence if $e_\Gamma\neq \gamma\in \Gamma$, then $L_\gamma(x)\neq x$ and the freeness of the action follows.  \qed
\vskip 0.1in
\noindent
Now we introduce our key notion: the {\bf landscape}. Landscapes are characterized by the height coordinates.
\begin{definition}
Let $y\in A^\Gamma$ be an element of totally finite height. We say that $y$ is a landscape if the following four conditions are satisfied.
\begin{itemize} 
\item If $d_G(\gamma,\delta)=1$ then $|H(y(\gamma))-H(y(\delta))|\leq 1$.
\item For all $n\geq 1$ there exists $M(y,n)>1$ such that if $H(y(\gamma))=n$ then there exists $\delta$, $d_G(\delta,\gamma)\leq M(y,n)$ so
that $H(y(\gamma))=1$.
\item for all $l\geq 1$ there exist $N(y,l)>1$ such that if $H(y(\gamma))=1$ then the ball $B_{N(y,l)}(G,\gamma)$ of radius $N(y,l)$ centered at $\gamma$ contains at least $l$ elements $\delta$ 
for which $H(y(\delta))=1$.
\item for all $m\geq 1$, there exists $S(y,m)>1$ so that every ball $B_{S(y,m)}(G,\delta)$ in our graph $G$ contains an element $\kappa$
such that $H(y(\kappa))\geq m$.
\end{itemize}
\end{definition}
\noindent
We call a landscape proper if its Cantor coordinates amount to
a proper Cantor labeling of $\Gamma$.  The following lemma is straightforward.
\begin{lemma} \label{regu} Landscapes are regular and
if $x$ is an element of totally finite height in the orbit closure of a (proper) landscape, then $x$ is a (proper) landscape with the same structure contants as $y$.
\end{lemma} 
\section{Landscapes and minimality}\label{propland}
The goal of this section is to prove the following proposition.
\begin{proposition} \label{legvege}
Let $y\in A^\Gamma$ be a proper landscape.
Then the orbit closure of $y$ contains an invariant set $Y_y\subset A^\Gamma$ homeomorphic to $\K$
such that the restricted action of $\Gamma$ on $Y_y$ is free and minimal.
\end{proposition}
\proof
For each pair of integers $m,n\geq 1$ we consider the finite set $CU^{m,n}_\Gamma$. An element $B$ of $CU^{m,n}_\Gamma$ is a labeling of 
the vertices of the
ball $B_m(G,e_\Gamma)$ by elements of the set
$$\{0,1\}^m\times\{l\in \N \mid |n-l|\leq m\}\,,$$ 
\noindent
in such a way that the second coordinate of the label of $e_\Gamma$ equals to $n$.
Let $x\in A^\Gamma$ be a proper landscape. For each $m\geq 1$, we have a map
$$\Theta^m_x:\Gamma\to \cup_{n=1}^\infty CU^{m,n}_\Gamma$$
constructed in the following way. First of all,
$\Theta^m_x(\gamma)$ will be an element of $CU_\Gamma^{m,H(x(\gamma))}$. 
Let $\delta \in B_m(G,e_\Gamma)$. Then, $\Theta^m_x(\gamma) (\delta)= (a,b)$, where
\begin{itemize}
\item $a= \left(C(x(\gamma\delta))\right)_m\,.$
\item $b=H(x(\gamma\delta))\,.$
\end{itemize}
\noindent
Let $\cal{B}$ denote the countable set $\cup_{m=1}^\infty\cup_{n=1}^\infty CU^{m,n}_\Gamma$. For each $x\in A^\Gamma$, we have a partition
$\cal{B}=I_x\cup J_x \cup K_x$, where
\begin{itemize}
\item $K_x$ is the subset of labeled balls $B$ in $\cal{B}$ such that
$\Theta^m_x(\gamma)\neq B$ if $\gamma\in\Gamma$. Here $m$ is the radius of $B$. 
\item $J_x$ is the subset of labeled balls $B$ in $\cal{B}$ such that
there exists  $K_B>0$ so that if $H(x(\gamma))=1$, then there exists $\delta$, $d_G(\delta,\gamma)\leq K_B$ for which
$\Theta^m_x(\delta)=B$. Again, $m$ is the radius of $B$.
\item $I_x$ is defined as $\cal{B}\backslash (K_x\cup J_x)$.
\end{itemize}
\noindent
The following lemma is trivial.
\begin{lemma} \label{csokken}
If $x$ is a proper landscape and $x'$ is an element of totally finite height in the orbit closure of $x$,
then
\begin{itemize}
\item $K_{x'}\supseteq K_x$.
\item $J_{x'}\supseteq J_x$.
\item $I_{x'}\subseteq I_x$.
\end{itemize}
\end{lemma}
\noindent 
\begin{lemma} \label{reduction}
Let $x\in A^\Gamma$ be a a proper landscape and $B\in I_x$ be a labeled ball of radius $m$.
Then there exists a proper landscape  $x'$ in the orbit closure of $x$ such that
$B\in K_{x'}$ and $H(x'(e_\Gamma))=1$. 
\end{lemma}
\proof
By the definition of $I_x$, we have a sequence $\{\gamma_k\}^\infty_{k=1}$ such that \\ $H(x(\gamma_k))=1$ and
if $d_G(\delta,\gamma_k)\leq k$ then $\Theta^m_x(\delta)\neq B$. Let $x'\in A^\Gamma$ be the limit point of the
sequence $\{L_{\gamma_k}(x)\}_{k=1}^\infty$ in $A^\Gamma$. Note that such limit point must exist by the landscape conditions and $x'$ is
again a proper landscape.
Then $B\in K_{x'}$ and $H(x'(e_\Gamma))=1$. \qed
\begin{lemma}
Let $y\in A^\Gamma$ be a proper landscape. Then we have an element $z\in A^\Gamma$ in the orbit closure of
$y$ such that
\begin{itemize}
\item $H(z(e_\Gamma))=1.$
\item The set $I_z$ is empty.
\end{itemize}
\noindent
(we will call such elements $z\in A^\Gamma$  {\bf minimal})
\end{lemma}
\proof
Let $I_y=\{B_1, B_2,\dots....\}$. Using Lemma \ref{csokken} and \ref{reduction}, we can inductively
construct a sequence $\{y_n\}^\infty_{n=1}$ in the orbit closure of $y$ such 
that 
\begin{itemize}
\item $H(y_n(e_\Gamma))=1$ and
\item $B_i\notin I_{y_n}$ if $1\leq i \leq n$.
\end{itemize}
\noindent
Let $z$ be a limit point of the sequence $\{y_n\}^\infty_{n=1}$.
Then $z$ is a proper landscape and the set $I_z$ is empty. \qed
\vskip 0.1in
\noindent
Now we can finish the proof of our proposition.
Let $z\in A^\Gamma$ be the minimal proper landscape in the previous lemma. The invariant subspace $Y_y$ is defined
as the set of elements of totally finite height in the orbit closure $\overline{O}(z)$. By Lemma \ref{regu}, all other elements of $\overline{O}(z)$
are of totally infinite height.
Let $t\in Y_y$.
It is enough to prove that $\overline{O}(t)$ contains $z$, that is,
for all $m\geq 1$ there exists
$y_m\in\overline{O}(t)$ such that $\Theta^m_z(e_\Gamma)=\Theta^m_{y_m}(e_\Gamma)\,.$
Let $\Theta^m_z(e_\Gamma)=B_m$. Since $B_m\in J_z$ there exists $K_m>0$
such that if $H(z(\delta))=1$, then
there exists $\gamma$ so that
\begin{itemize}
\item $d_G(\delta,\gamma)\leq K_m$ and
\item $\Theta^m_z(\gamma)=B_m$.
\end{itemize}
Since $t\in \overline{O}(z)$, if $H(t(\delta))=1$ then
there exists some $\rho_m\in\Gamma$ such that
$d_G(\rho_m,\delta)\leq K_B$ and $\Theta^m_t(\rho_m)=B_m$. That is,
$\Theta^m_{y_m}(e_\Gamma)=B_m$ if $y_m=L_{\rho_m}(t)\,.$

\noindent
Finally, we need to show that $t$ is not an isolated point in $Y_y$.
Let $\Theta^m_t(e_\Gamma)=B'_m.$ By our third landscape condition and the minimality of $z$, there exists
$e_\Gamma\neq \gamma_m\in\Gamma$ such that
$\Theta^m_t(\gamma_m)=B'_m$ as well. Hence, $\Theta^m_{L_{\gamma_m}(t)}=B'_m\,.$
Since by freeness $ L_{\gamma_m}(t)\neq t$ for all $m\geq 1$, we have that
$t$ is not an isolated point. \qed

\section{Hilly landscapes and Borel hyperfiniteness}
Let $z\in A^\Gamma$ be a proper landscape. We say that $z$ is {\bf hilly} if for all $n\geq 1$ there exists
$Q_n$ such that the induced
graph in $G$ on the set $W^n_z\subset \Gamma$, where
$$W^n_z=\{\gamma\,\mid \, H(z(\gamma))\leq n\}$$
\noindent
has components of size at most $Q_n$.
Clearly, if $y$ is a minimal landscape in the orbit closure of a hilly landscape $z$ then
$y$ is hilly as well (with the same structure constants $\{Q_n\}^\infty_{n=1}$). Let $\alpha:\Gamma\actson X$ be a Borel action of $\Gamma$
on a Borel space $X$. We say that $p,q\in X$ are equivalent, $p\equiv_E q$, if for some $\gamma\in\Gamma$, $\alpha(\gamma)(p)=q.$ 
The equivalence relation $E$ is called the orbit equivalence relation of the action. Recall that $\alpha$ is called Borel hyperfinite, if $E$
is the increasing union of some finite Borel equivalence relations $E_1\subset E_2\subset\dots$.
\begin{proposition}
Let $y$ be a minimal hilly landscape. Then the action of $\Gamma$ on the totally finite part $Y$ of the orbit closure of $y$ is 
Borel hyperfinite (consequently amenable, see Section \ref{code}).
\end{proposition}
\proof
We define the finite equivalence relation $E_n$ on $Y$ in the following way.
If $t,s\in Y$ then $t\equiv_{E_n}s$ if
\begin{itemize}
\item either $t=s$,
\item or $L_\gamma(t)=s$ for some $e_\Gamma\neq \gamma\in\Gamma$, such that there exists a path $(e_\Gamma=\gamma_1,\gamma_2,\dots,\gamma_l=\gamma)$ in $G$
for which $t(\gamma_i)\leq n$ holds if $1\leq i \leq l$.
\end{itemize}
\noindent
Since the elements of $Y$ are hilly, $E_n$ is indeed a finite Borel equivalence relation.
Clearly, $E_1\subset E_2\subset\dots$ and $\cup^\infty_{n=1} E_n$ is the orbit equivalence relation 
on $Y$. Therefore the action of $\Gamma$ on $Y$ is Borel hyperfinite.
\qed
\begin{proposition}
There are hilly landscapes on $\Gamma$.
\end{proposition}
\proof
We use a fractal-like construction to build the landscape (one should note that the so-called $(C,F)$- 
construction in \cite{Dan1} also has a fractal-like character). So, let $A_0=\{e_\Gamma\}$, $A_1=A_0\cup \{\gamma_1\}$,
where $d_G(e_\Gamma,\gamma_1)=30$. Let $A_2=A_1\cup \gamma_2 A_1$, where
$d_G(e_\Gamma, \gamma_2)=300$ and inductively, let $A_n=A_{n-1}\cup\gamma_n A_{n-1}$, where
$d_G(e_\Gamma,\gamma_n)=3 (10^n).$ Let $A=\cup_{n=1}^\infty A_n$. Observe that
any non-unit element of $A$ can be uniquely written as $\gamma_{n_k}\gamma_{n_{k-1}}\dots \gamma_{n_1}$, where
$n_k>n_{k-1}>\dots >n_1$. We define the subsets
$$A\supset Q_1 \supset Q_2\supset \dots$$
\noindent
by
$$Q_n:= \{e_\Gamma\}\cup \{\delta\,\mid\, \delta=\gamma_{n_k}\gamma_{n_{k-1}}\dots \gamma_{n_1} \,\mbox{and}\, n_1\geq n\}\,.$$
\noindent
So, in particular $\gamma_n\in Q_n$.  Observe that
\begin{itemize}
\item If $\gamma,\delta\in Q_n$ then $B_{10^n}(G,\delta)$ and $B_{10^n}(G,\gamma)$ are disjoint.
\item If $\gamma\in Q_k$ and $\delta\in Q_l$, where $k<l$ then either
$B_{10^k}(G,\gamma)\subset B_{10^l}(G,\delta)$ or
$B_{10^k}(G,\gamma)\cap B_{10^l}(G,\delta)=\emptyset$.
\end{itemize}
\noindent
Now we can define the landscape $z$ on $\Gamma$ in the following way. Let $H(z(\gamma))=l$ if
\begin{itemize}
\item $\gamma\in B_{10^l}(\delta)$ for some $\delta\in Q_l$ and also
\item $\gamma\notin B_{10^k}(\rho)$ if $\rho \in Q_k$ and $k<l$.
\end{itemize}
\noindent
Also, let $\gamma\to C(z(\gamma))$ be an arbitrary proper labeling.
It is easy to check that $z$ is, in fact, a hilly landscape. \qed
\vskip 0.1in
\noindent
\begin{remark} We can construct an explicit hilly landscape on the group of integers $Z$ using the construction above.
Call a non-negative integer $n$ {\it ternary} if all the digits of $n$ are $0$ or $3$:  $0,3,30,33,300,303,\dots$ We only need to define $H:\Z\to \N$.
Let $H(n)=1$ if $n$ is a ternary number. In general, let  $H(n)=k$, if $k$ is the smallest non-negative integer such that $|n-t|\leq 10^k$, where $t$ is a ternary
number and $10^k$ divides $t$.
\end{remark}

\section{Landscapes with rivers} \label{river}
The Borel hyperfinite construction of the previous section cannot be extended
to a purely infinite action, so we need a different idea.
Let $G$ be the Cayley graph of $\Gamma$ as in the previous sections. Also, let $T$ be the infinite tree for which
all the vertex degrees are four (the $4$-tree). 
A {\bf river}
is a bilipschitz embedding of the $4$-tree $T$ into $G$, that is, a map $\Psi:V(T)\to\Gamma$ such
that there exists some $C>0$ so that for all $x,y\in V(T)$
$$C^{-1}d_T(x,y)\leq d_G(\Psi(x),\Psi(y))\leq C d_T(x,y)\,.$$
\noindent
Let us also assume that $H(y(\gamma)):=d_G(\Psi(V(T),\gamma))+1$ defines a landscape on $\Gamma$ and $e_\Gamma\in \Psi(V(T))$.
We call such $y$ a {\bf landscape with river}. 
\begin{proposition} \label{riverexist}
Landscapes with rivers do exist on non-amenable groups. 
\end{proposition}
\proof
By Theorem 1.5 of \cite{BS}, bilipschitz embeddings $\Psi_1:V(T)\to\Gamma$ exist. Clearly,
the resulting element $y$ would satisfy the first three landscape conditions. However, the fourth condition might not be satisfied, say, because $\Gamma$ is the
free group and $\Psi_1$ is surjective.
So, let us consider a bilipschitz map $\Phi:V(T)\to V(T)$ such that
for all $t\in V(T)$, there is at least one branch $B_T$ of $t$ in the tree $T$ so that
$B_T\cap \Phi(V(T))$ is empty. Let $\Psi=\Psi_1\circ\Phi$. We can also assume that $e_\Gamma\in \Psi(V(T))$.
 Then the resulting element $y$ will
satisfy the fourth landscape condition as well. \qed
\vskip 0.1in
\noindent
The following proposition will be crucial in the next section. Note that for a set $A$, $\mbox{Fin}(A)$ denotes the
the family of all finite subsets of $A$.
\begin{proposition} \label{PropertyA} 
Let $\Psi:V(T)\to \Gamma$ be a bilipschitz embedding of the $4$- tree into our Cayley graph $G$ in such a way that $e_\Gamma\in \Psi(V(T))$ and \\
$H(y(\gamma)):=d_G\left(\Psi(V(T)),\gamma\right)+1$ defines
a landscape on $\Gamma$.
Let $C>0$ be an integer such that if $x,y\in V(T)$ then
$$C^{-1} d_T(x,y)\leq d_G(\Psi(x),\Psi(y))\leq C d_T(x,y)\,.$$
\noindent
Then for all $m\geq 1$ we have a map
$\kappa_m:\Gamma\to\mbox{Fin}\left(\Psi(V(T))\right)$ such that 
\begin{itemize}
\item For all $\gamma\in\Gamma$, $|\kappa_m(\gamma)|=m.$
\item For all $\gamma\in\Gamma$, $\kappa_m(\gamma)\subset B_{d_G(\Psi(V(T)),\gamma)+Cm}(G,\gamma)\,.$
\item
If $d_G(\gamma_1,\gamma_2)=1$, then
\begin{equation}\label{egyenlet}
|\kappa_m(\gamma_1)\triangle \kappa_m(\gamma_2)|\leq 2(d_G(\Psi(V(T)),\gamma_1)+1)C\,.
\end{equation}
\end{itemize}
\end{proposition}
\proof
We use the classical construction that shows the $4$-tree has Property A. This process hopefully explains why we call these objects rivers.
First, let us fix an infinite ray $\{t_i\}^\infty_{i=0}$ in $V(T)$. That is, 
$d_T(t_0,t_i)=i$ and $d_T(t_{i-1},t_i)=1\,.$  If $s\in V(T)$, then we have a unique path $(s_1,s_2,\dots,s_l)$ such
that $s=s_1,s_l=t_i$ for some non-negative integer $i$, and $s_{l-1}\notin \{t_i\}^\infty_{i=0}$.
Then, we consider the infinite path $P(s)=(s_1,s_2,\dots,s_l,s_{l+1}\dots)$, where for all $j\geq 1$, $s_{l+j}=t_{i+j}$.
So, for each $m\geq 1$, we have the path $P_m(s)=(s_1,s_2,\dots,s_m)$.
Then for all $s\in V(T)$, $|P_m(s)|=m$. Also, if $p,q\in V(T)$ and $d_T(p,q)=a$, then $|P_m(p)\triangle P_m(q)|\leq a$.
Now for all $\gamma\in\Gamma$, pick an element $\delta_\gamma$ in $\Psi(V(T))$ such 
that $d_G(\delta_\gamma,\gamma)=d_G(\Psi(V(T)),\gamma)$ and let
$$\kappa_m(\gamma)=\Psi\left(P_m(\Psi^{-1}(\delta_\gamma))\right)\,.$$
Now, if $d_G(\gamma_1,\gamma_2)=1$, then $d_G(\delta_{\gamma_1},\delta_{\gamma_2})\leq 2(d_G(\Psi(V(T)),\gamma_1)+1)\,.$
Hence, $|\kappa_m(\gamma_1) \triangle \kappa_m(\gamma_2)|\leq 2(d_G(\Psi(V(T)),\gamma_1)+1)C.$
Also, for all $\gamma\in\Gamma$, $\kappa_m(\gamma)\subset B_{d_G(\Psi(V(T)),\gamma)+Cm}(G,\gamma)\,.$ \qed
\vskip 0.1in
\noindent
Let $y\in A^\Gamma$ be an element of totally finite height such that
$H(y(\gamma))=d_G(\Psi(V(T)),\gamma)+1$ for some river with bilipschitz constant $C$.
Let 
$$H_1(y)=\Psi(V(T))=\{\gamma\in\Gamma\mid H(y(\gamma))=1\}.$$
\noindent
Let $G_y$ be a graph on the vertex set $H_1(y)$ defined in the following way.
\begin{itemize}
\item $V(G_y)=H_1(y)\,.$
\item $(p,q)\in E(G_y)$ if $d_G(p,q)\leq C$.
\end{itemize}
\noindent
Then $G_y$ has bounded vertex degrees and it is quasi-isometric to the $4$-tree $T$. In particular, $G_y$ has positive
Cheeger constant.
Recall that the Cheeger constant of an infinite graph $J$ is defined in the following way.
$$c(J)=\inf_{H\in Fin(V(J))} \frac{|\partial (H)|}{|H|}\,,$$
where $\partial(H)=\{p\in H\mid \exists q\notin H, d_J(p,q)=1\}.$
Now let $z\in A^\Gamma$ be an element of totally finite height in the orbit closure of $y$. We can construct
the graph $G_z$ on $H_1(z)=\{\gamma\in\Gamma\mid H(z(\gamma))=1\}$ as above.
\begin{lemma} \label{cheeger}
The graph $G_z$ is connected and is of bounded vertex degree. Also, $G_z$ has positive Cheeger constant, in fact, $c(G_z)\geq c(G_y)\,.$ 
\end{lemma}
\proof Since $K_z\subseteq K_{z'}$ it is clear that $G_z$ has bounded vertex degrees and  $c(G_z)\geq c(G_y)$.  Now we prove 
the connectivity of $G_z$.
Let $\gamma,\gamma\delta\in H_1(z)$ such that $\delta\in B_m(G,e_\Gamma)$. Note that
if $p,q\in H_1(y)$ and $d_G(p,q)\leq m$, then $d_{G_y}(p,q)\leq Cm$.  Since $z$ is in the orbit closure of $y$, there exist 
$\rho,\rho\delta\in H_1(y)$
such that
$$\Theta^{C^2 m}_z(\gamma)=\Theta^{C^2 m}_y(\rho)\,.$$
\noindent
Since $\rho$ and $\rho\delta$ can be connected by a path in $H_1(y)$ inside the ball $B_{C^2 m}(G,\rho)$ we can conclude that
$\gamma$ and $\gamma\delta$ can be connected by a path in $H_1(z)$ inside the ball $B_{C^2 m}(G,\gamma)$. This finishes
the proof our lemma. \qed

\section{The Cantor code for amenability} \label{code}
Let $y$ be a proper landscape with river (so we also assume
that a proper labeling $\gamma\to C(y(\gamma))$ is given) and for each 
$m\geq 1$ let
$\kappa_m:\Gamma\to \mbox{Fin}(\Psi(V(T)))$ be the map as in Section \ref{river}.
That is, for all $\gamma\in\Gamma$
$$\kappa_m(\gamma)\subset B_{F_{m,n}}(G,\gamma)\,,$$
where $n=H(y(\gamma))=d_G(\Psi(V(T)),\gamma)$ and $F_{m,n}=Cm+n$.
For $\gamma\in\Gamma$, let $L_m(\gamma)\subset B_{F_{m,n}}(G,e_\gamma)$
be the subset such that
$$\gamma L_m(\gamma)=\kappa_m(\gamma)\,.$$
For each $m,n\geq 1$ 
let $\{a^{m,n}_1,a^{m,n}_2,\dots,a^{m,n}_{\tau_{m,n}}\}$ be an enumeration of the 
finite subsets of
$B_{F_{m,n}}(G,e_\Gamma)$, where $\tau_{m,n}=2^{|B_{F_{m,n}}(G,e_\Gamma)|}\,.$
So, for each $\gamma\in\Gamma$ we have an element
$c_\gamma\in\{0,1\}^\N$ constructed in the following way.
Let $L_m(\gamma)=a^{m,n}_{i_{m,n,\gamma}}$, where
$1\leq i_{m,n,\gamma}\leq \tau_{m,n}$.
Now let $c_{m,\gamma}$ be the concatenation of $i_{m,n,\gamma}$ pieces of 
the string $010$.
Let
$$c_\gamma=(11 c_{1,\gamma} 11 c_{2,\gamma} 11 c_{3,\gamma} 11\dots)\in \{0,1\}^N\,.$$
\noindent
Therefore, for each $\gamma\in\Gamma$ we have two elements
of the Cantor set:
$$C(y(\gamma))=(u^\gamma_1  u^\gamma_2 \dots)$$
\noindent
and
$$c_\gamma=(v^\gamma_1 v^\gamma_2 \dots)\,.$$
Let $z\in A^\Gamma$ be defined by
\begin{itemize}
\item $H(z(\gamma))=H(y(\gamma))\,.$
\item $C(z(\gamma))=(u^\gamma_1v^\gamma_1u^\gamma_2 v^\gamma_2u^\gamma_3v^\gamma_3\dots)$
\end{itemize}
\noindent
Clearly,
$z$ is a proper landscape. Notice that $z$ encodes the landscape $y$ and for each $m\geq 1$, the
system $\{\kappa_m(\gamma)\}_{\gamma\in\Gamma}$. Finally, let $x$ be a minimal element in the orbit closure of $z$ and $Y$ be
the totally finite part of the orbit closure of $x$.
\begin{proposition} \label{freeminame}
The action of $\Gamma$ on $Y$ is free, minimal and amenable.
\end{proposition}
\proof By Proposition \ref{legvege}, freeness and minimality follow.
So, let us recall the definition of amenable actions on  locally compact spaces.
\begin{definition} \cite{Anan} Let $\Gamma$ be a finitely generated group
with a finite generating system $\Sigma$. Let $\alpha:\Gamma \actson X$ be a continuous action of $\Gamma$
on the locally compact space $X$. The action $\alpha$ is topologically amenable if there exists
a sequence $\{g_m:X\times\Gamma\to\R\}^\infty_{m=1}$ of non-negative functions such that
\begin{itemize}
\item For all $m\geq 1$ and $p\in X$, $\sum_{\gamma\in\Gamma} g_m(p,\gamma)=1.$
\item for all generator $\sigma\in \Sigma$
$$\sum_{\gamma\in\Gamma} |g_m(\alpha(\sigma)(p),\sigma\gamma)-g_m(p,\gamma)|$$
\noindent
uniformly tends to zero on the compact subsets of $X$.
\end{itemize}
\end{definition}
\noindent
Let $t\in Y$. Then 
$$C(t(\gamma))=(u^1_{t,\gamma} v^1_{t,\gamma}u^2_{t,\gamma} v^2_{t,\gamma}\dots)\,.$$
\noindent
Let us consider
$$C_v(t(\gamma))=(v^1_{t,\gamma}, v^2_{t,\gamma}\dots)\,.$$
\noindent
By our construction,
$$C_v(t(\gamma))=(11 d^1_{t,\gamma} 11 d^2_{t,\gamma} 11\dots)\,,$$
\noindent
where
$d^m_{t,\gamma}$ is the concatenation of $j_{m,t,\gamma}$ pieces of the string $010$.
Also, $j_{m,t,\gamma}\leq \tau_{m,H(t(\gamma))}\,.$
Let us define $g_m:Y\times \Gamma\to\R$ in the following way.
Let $g_m(t,\rho)=\frac{1}{m}$ if
$\rho\in a^{m,H(t(e_\Gamma))}_{j_{m,t,e_\Gamma}}$,
otherwise, let $g_m(t,\rho)=0$.
Clearly, $g_m$ is continuous and for all $t\in Y$, $\sum_{\rho\in\Gamma} g_m(t,\rho)=1\,.$
Since $Y$ is contained in the orbit closure of the element $z$, for all $t\in Y$ there exists $\delta\in\Gamma$ 
such that
\begin{itemize}
\item
$H(z(\delta))=H(t(e_\Gamma))\,.$
\item
$L_m(\delta)=a^{m,H(t(e_\Gamma))}_{j_{m,t,e_\Gamma}}\,.$
\item
$L_m(\delta\sigma)= a^{m,H(L_\sigma(t)(e_\Gamma))}_{j_{m,L_\sigma(t),e_\Gamma}}\,.$
\end{itemize}
\noindent
Therefore by \eqref{egyenlet},
$$\sum_{\rho\in\Gamma}|g_m(L_\sigma(t),\sigma\rho)-g_m(t,\rho)|\leq \frac{2(H(t(e_\Gamma))+2)C}{m}\,.$$
That is,
$$\sum_{\rho\in\Gamma}|g_m(L_\sigma(t),\sigma\rho)-g_m(t,\rho)|$$
uniformly tends to zero on the set $Y_n=\{t\in Y\mid H(t(e_\gamma))\leq n\}$.
Since for all compact open set $K\subset Y$ there exists $n\geq 1$ such that $K\subset Y_n$,
our proposition follows. \qed

\section{The combinatorial version of paradoxicality}
Let $z\in A^\Gamma$ be a minimal proper landscape and for $m\geq 1$ let $\Theta^m_z:\Gamma\to \cup_{n=1}^\infty CU^{m,n}_\Gamma$ be the map
defined in Section \ref{propland}. We say that the subset $T\subseteq \Gamma$ is $z$-local if there
exists $m\geq 1$ and a finite subset $S\subset \cup_{n=1}^\infty CU^{m,n}_\Gamma$ such that
$T=(\Theta^m_z)^{-1}(S)\,.$ Notice that the $z$-locality of the subset $T$
means that the membership of $T$ can be locally verified.
We call a $z$-local subset $T$ $z$- paradoxical if there
exist pairwise disjoint $z$-local subsets $T_1, T_2,\dots, T_{p+q}$ and elements
$\gamma_1,\gamma_2,\dots,\gamma_{p+q}\in\Gamma$ such that 
$T_j\subset T$ for all $j$, and
$$T=\cup_{i=1}^p T_i\gamma_i= \cup_{j=p+1}^{p+q} T_j\gamma_j\,.$$
\noindent
Now let $Z$ be the totally finite height part of the orbit closure $\overline{O}(z)$ of $z$.
We define the map $\Theta^m_Z:Z\to \cup_{n=1}^\infty CU^{m,n}_\Gamma$ by
$$\Theta^m_Z(x)=\Theta^m_x(e_\Gamma)\,.$$
\noindent
Note that $\Theta^m_Z$ is a locally constant function, hence if $S\subset \cup_{n=1}^\infty CU^{m,n}_\Gamma$ is a finite subset then
$(\Theta^m_{Z})^{-1}(S)$ is a compact-open subset of the locally compact space $Z$.
Moreover, by the definition of the product topology, any compact-open subset $U$ of $Z$ can be written
as $(\Theta^m_{Z})^{-1}(S)$ for some $m\geq 1$ and finite subset $S\subset \cup_{n=1}^\infty CU^{m,n}_\Gamma$.
The key observations of this section are the following propositions.
\begin{proposition}\label{key}
Let $m\geq 1$ and let $S\subset \cup_{n=1}^\infty CU^{m,n}_\Gamma$ be a finite subset. Let $z\in A^\Gamma$ be a minimal proper landscape and
$Z$ be as above. Suppose that the $z$-local subset $T=(\Theta^m_z)^{-1}(S)$ is $z$-paradoxical.
Then $U=(\Theta^m_Z)^{-1}(S)$ is a paradoxical compact-open subset of $Z$. Consequently, if all $z$-local subsets of $\Gamma$
are $z$-paradoxical then the action of $\Gamma$ on $Z$ is purely infinite.
\end{proposition}
\proof
Let $\gamma_1,\gamma_2,\dots,\gamma_{n+m}\in\Gamma$ such that $T_j\subset T$ for all $j$ and
$$T=\cup_{i=1}^p T_i\gamma_i= \cup_{j=p+1}^{p+q} T_j\gamma_j\,.$$
\noindent
Then there exists $l>m$ and $S_1, S_2,\dots, S_{n+m}\in\mbox{Fin}
(\cup_{n=1}^\infty CU^{l,n}_\Gamma)$ such that
for all $1\leq i \leq p+q$, $T_i=(\Theta^l_z)^{-1}(S_i)\,.$ Now observe that
$$U=\cup_{i=1}^p L_{\gamma_i} (U_i)= \cup_{j=p+1}^{p+q} L_{\gamma_j} (U_j)\,,$$
\noindent
where $U_i=(\Theta^l_Z)^{-1}(S_i)\,.$ Hence, $U$ is indeed paradoxical. \qed
\begin{proposition}
\label{clean}
Let $z\in A^\Gamma$ be a landscape and let $w\in A^\Gamma$ be an element
of totally finite height in the orbit closure of $z$. Let $m>0$ and
let $S\in\mbox{Fin}(\cup^\infty_{n=1} CU^{m,n}_\Gamma)$. Also, let $l>m$, for $1\leq i \leq p+q$ 
let $S_i\in \mbox{Fin}(\cup^\infty_{n=1}CU^{l,n}_\Gamma)$ and $\gamma_1,\gamma_2,\dots,\gamma_{p+q}\in\Gamma$
such that
\begin{itemize}
\item
$T=(\Theta^m_z)^{-1}(S), T_i=(\Theta^l_z)^{-1}(S_i)$.
\item
$T_i\subset T$ for $1\leq i \leq p+q$.
\item
$T=\cup^p_{i=1} T_i\gamma_i=\cup^{p+q}_{j=p+1} T_j\gamma_j$.
\end{itemize}
\noindent
Then the sets $\{T^w_i\}^{p+q}_i$ are disjoint, $T^w_i\subset T^w$ and $T^w=\cup^p_{i=1} T^w_i\gamma_i=\\ \cup^{p+q}_{j=p+1} T^w_j\gamma_j$,
where $T^w=(\Theta^m_w)^{-1}(S), T_i=(\Theta^l_w)^{-1}(S_i)$. That is,
$T^w$ is $w$-paradoxical (note that empty sets are paradoxical by definition).
\end{proposition}
\proof
Let $a$ be an integer such that
$\gamma_1,\gamma_2,\dots,\gamma_{n+m}\in B_a(G,e_\Gamma)$. First, let us prove that
$T^w_i\cap T^w_j=\emptyset$ if $i\neq j$.
Suppose that $\gamma\in T^w_i\cap T^w_j$. Since $w$ is in the orbit closure
of $z$, there exists $\delta\in\Gamma$ such that $\Theta^l_z(\delta)=\Theta^l_w(\gamma)$.
Hence, $\delta\in T_i\cap T_j$ leading to a contradiction.
\noindent
Now let $\gamma\in T_i^w$. We need to show that
$\gamma\gamma_i\in T^w$. Again, we have $\delta\in\Gamma$ such that
$\Theta^{a+l+m}_z(\delta)=\Theta^{l+a+m}_w(\gamma)\,.$
Then $\delta\in T_i$, so $\delta\gamma_i\in T$, hence $\gamma\gamma_i\in T^w$.
\noindent
Finally, let $\gamma\in T^w$. Let us show that there exists
$1\leq i \leq p$ such that $\gamma\gamma_i^{-1}\in T^w_i$ and
$p+1\leq j \leq p+q$ such that $\gamma\gamma^{-1}_j\in T^w_j$.
Again, let $\delta\in\Gamma$ such that
$\Theta^{a+l+m}_z(\delta)=\Theta^{l+a+m}_w(\gamma)\,.$
Then $\delta\in T$, hence for some $1\leq i \leq p$ and
$p+1\leq j \leq p+q$ we have that $\delta\gamma_i^{-1}\in T_i$ and
$\delta\gamma_j^{-1}\in T_j$. Thus,
$\gamma\gamma_i^{-1}\in T^w_i$ and
 $\gamma\gamma^{-1}_j\in T^w_j$. \qed

\section{Paradoxicalization}
Let $z,z'\in A^\Gamma$ be landscapes and let $l\geq 1$ be an integer.
Then $z\equiv_l z'$ if for all $\gamma\in\Gamma$:
\begin{itemize}
\item $C(z(\gamma))_l=C(z'(\gamma))_l$.
\item
If $C(z(\gamma))=(a_1b_1a_2b_2\dots)$ and $C(z'(\gamma))=(c_1d_1c_2d_2\dots)$, then
for all $n\geq 1$, $a_n=c_n$. 
\item $H(z(\gamma))=H(z'(\gamma))$.
\end{itemize}
\noindent
The following lemma is a straightforward consequence of the definition.
\begin{lemma} \label{modi}
Let $z\in A^\Gamma$ be a landscape, let $m>0$ and $S\in \mbox{Fin}(\cup_{n=1}^\infty CU^{m,n}_\Gamma)$. Also, let $l>m$ and 
for $1\leq i \leq p+q$, let $S_i\in \mbox{Fin}(\cup_{n=1}^\infty CU^{l,n}_\Gamma)$ and $\gamma_1, \gamma_2\dots, \gamma_{p+q} \in\Gamma$ 
such that
$\emptyset \neq T$ is $z$-paradoxical and
$$ T=(\Theta^m_z)^{-1}(S)= \cup_{i=1}^p T_i\gamma_i=\cup_{j=p+1}^{p+q} T_j\gamma_j\,,$$
where $T_i=(\Theta^l_z)^{-1}(S_i)$. Then if $z\equiv_r z'$, where $r\geq l$ : $T$ is $z'$-paradoxical as well.
\end{lemma}
\noindent 
One of main ingredients of the proof of Theorem \ref{fotetel}
is the following proposition.
\begin{proposition} \label{ingr}
Let $z$ be a minimal landscape and let $\emptyset \neq T=(\Theta^m_z)^{-1}(S)$ be a $z$-local set,
where $S\in \mbox{Fin}(\cup_{n=1}^\infty CU^{m,n}_\Gamma)$. Let $m\leq m'$. Then there exists
$z'\in A^\Gamma$ such that $z'\equiv_{m'} z$ and $T$ is $z'$-paradoxical. 
\end{proposition}
\proof
Since $z$ is minimal there exists some $R_T>1$ such that if $\gamma\in\Gamma$ then
$B_{R_T}(G,\gamma)\cap T\neq \emptyset$.
Let us construct a graph $G_T$ with vertex set $T$ in the following way.
The vertices $p,q\in T$ are adjacent in $G_T$ if and only if
$d_G(p,q)\leq 3R_T$.
It is easy to see that $G_T$ is a connected graph with bounded vertex degrees and $G_T$ is quasi-isometric
to $G$. Since $G$ is the Cayley graph of a non-amenable group, $G$ has positive Cheeger constant. So, since $G$ and
$G_T$ are quasi-isometric, $G_T$ has positive Cheeger constant as well (Theorem 18.13 \cite{Dru}).
Therefore, by the main result of \cite{Deuber} $G_T$ is a paradoxical graph. That is, there exist injective maps
$\phi:T\to T$, $\psi:T\to T$ and $K>0$ such that 
\begin{itemize}
\item $\phi(T)\cap\psi(T)=\emptyset.$
\item For every $x\in T$, $d_G(x,\phi(x))<K$, $d_G(x,\psi(x))<K$.
\end{itemize}
\noindent
Therefore, there exist elements $\gamma_1,\gamma_2,\dots, \gamma_{p+q}\in\Gamma$ such that for
any $x\in T$ 
\begin{itemize}
\item
there exists $1\leq i \leq p$ such that $\phi(x)\gamma_i=x$,
\item
there exists $p+1\leq j \leq q$ such that $\psi(x)\gamma_j=x$.
\end{itemize}
\noindent
For $1\leq i \leq p$ let
$$T_i=\{y\in T\mid \,\mbox{there exists $x\in T$ such that $\phi(x)\gamma_i=x$}\}\,.$$
\noindent
For $p+1\leq j  \leq p+q$ let
$$T_j=\{y\in T\mid \,\mbox{there exists $x\in T$ such that $\psi(x)\gamma_j=x$}\}\,.$$
\noindent
Then $T_1,T_2,\dots T_{p+q}$ are disjoint sets and
$$T=\cup^p_{i=1} T_i\gamma_i=\cup^{p+q}_{p+1} T_j\gamma_j\,.$$
\noindent
We need to construct $z'\in A^\Gamma$ such that $z\equiv_{m'} z'$ for all 
$1\leq i \leq p+q$, $T_i$ is $z'$-local.
Let $a_1< a_2 <\dots < a_{n+m}$ be consecutive even numbers such that $m'<a_1$. 
For $\gamma\in\Gamma$ let $H(z'(\gamma))=H(z(\gamma))$ and
\begin{itemize}
\item If $\gamma \notin \cup_{i=1}^{p+q} T_i$, then for all $1\leq i \leq p+q$ let the $a_i$-th Cantor coordinate
of $z'(\gamma)$ be $0$.
\item If $\gamma\in T_j$ then let the $a_j$-th Cantor coordinate of $z'(\gamma)$ be $1$ and if $i\neq j$ then
let the $a_i$-th Cantor coordinate
of $z'(\gamma)$ be $0$.
\end{itemize}
\noindent
For $l>a_{p+q}$ the $l$-th Cantor coordinates of the elements of $\Gamma$ will be chosen in such a way to
make $z'$ proper. For $l< a_1$ let the $l$-th Cantor coordinate of $z'(\gamma)$ be equal to the
$l$-th Cantor coordinate of $z(\gamma)$. It is easy to see that for the resulting proper landscape $z'$, 
$z\equiv_{m'} z$ and for all $1\leq i \leq p+q$, $T_i$ is $z'$-local.
\qed

\section{The proof of Theorem \ref{fotetel}}
Now we are in the position to prove our theorem.
First, let $\{S_i\}^\infty_{i=1}$ be an enumeration of the set $\cup^\infty_{m=1} \mbox{Fin}(\cup^\infty_{n=1} CU^{m,n}_\Gamma)$ and
$S_i\in\mbox{Fin}(\cup^\infty_{n=1} CU^{m_i,n}_\Gamma)$.
Let $y$ be a minimal landscape with a river. Let $z$ and $x$ be the landscapes as in Proposition \ref{freeminame}. Finally, we define $t_0\in A^\Gamma$ in
the following way. 
\begin{itemize}
\item For all $\gamma\in\Gamma$, $H(t_0(\gamma))=H(x(\gamma))$. 
\item For $\gamma\in\Gamma$, if $C(x(\gamma))=(u_1u_2u_3\dots)$,
let $C(t(\gamma))=(u_10u_20 u_3 0\dots)$.
\end{itemize}
We define $\pi_{odd}:A^\Gamma\to A^\Gamma$ as follows.
For $p\in A^\Gamma$ and $\gamma\in\Gamma$
\begin{itemize}
\item $H(\pi_{odd}(p)(\gamma))=H(p(\gamma))\,.$
\item If $C(p(\gamma))=(v_1v_2v_3v_4v_5\dots)$ then
$C(\pi_{odd}(p)(\gamma))=(v_1v_3v_5\dots)$.
\end{itemize}
\noindent
So, $\pi_{odd}(t_0)=x$ and $t_0$ is again a minimal landscape.
Now we start our inductional process.

\noindent
{\bf Step 0.}
If $(\Theta^{m_1}_{t_0})^{-1} (S_1)$ is empty then let $l_1=m_1, p_1=0, q_1=1, S^1_1=\emptyset, \\ \gamma_1=e_\Gamma, t_1=t_0$.
If $(\Theta^{m_1}_{t_0})^{-1} (S_1)$ is a (non-empty) $t_0$-local set, then
let $t'_0\in A^\Gamma$ be such that $t'_0\equiv_{m_1} t_0$  and $S_1$ is $t'_0$-paradoxical (Proposition \ref{ingr}). Also, 
let $\gamma^1_1,\gamma^1_2,\dots,\gamma^1_{p_1+q_1}\in\Gamma$
and $S^1_1,S^1_2,\dots,S^1_{p_1+q_1}\in \mbox{Fin}( \cup^\infty_{n=1} CU^{l_1,n})$  such that 
\begin{itemize}
\item For all $1\leq i \leq p_1+q_1$, $(\Theta^{l_1}_{t'_0})^{-1}(S^1_i)\subset (\Theta^{m_1}_{t'_0})^{-1}(S_1)$.
\item The sets $\{(\Theta^{l_1}_{t'_0})^{-1}(S^1_i)\}_{i=1}^{p_1+q_1}$ are disjoint.
\item $(\Theta^{m_1}_{t'_0})^{-1}(S_1)= \cup_{i=1}^{p_1} ((\Theta^{l_1}_{t'_0})^{-1}(S^1_i))\gamma^1_i=\cup_{j=p_1+1}^{p_1+q_1}( (\Theta^{l_1}_{t'_0})^{-1}(S^1_j))\gamma^1_j \,.$
\end{itemize}
\noindent
Finally, let $t_1$ be a minimal landscape in the orbit closure of $t'_0$.
Then, 
\begin{itemize}
\item For all $1\leq i \leq p_1+q_1$, $(\Theta^{l_1}_{t_1})^{-1}(S^1_i)\subset (\Theta^{m_1}_{t_1})^{-1}(S^1_1)$.
\item The sets $\{(\Theta^{l_1}_{t_1})^{-1}(S^1_i)\}_{i=1}^{p_1+q_1}$ are disjoint.
\item $(\Theta^{m_1}_{t_1})^{-1}(S^1_1)= \cup_{i=1}^{p_1} ((\Theta^{l_1}_{t_1})^{-1}(S^1_i))\gamma^1_i=\cup_{j=p_1+1}^{p_1+q_1} ((\Theta^{l_1}_{t_1})^{-1}(S^1_j))\gamma^1_j \,.$
\end{itemize}
\noindent
So, 
\begin{itemize}
\item $t_1$ is a minimal landscape.
\item $\pi_{odd}(t_1)$ is in the orbit closure of $x$.
\item $(\Theta^{m_1}_{t_1})^{-1}(S_1)$ is $t_1$-paradoxical.
\end{itemize}
\vskip 0.1in
{\bf Step k.}
Suppose that we have a minimal landscape $t_k$ and we also have
\begin{itemize}
\item For all $1\leq a \leq k$,   $\gamma^a_1,\gamma^a_2,\dots,\gamma^a_{p_a+q_a}\in\Gamma$. 
\item For all $1\leq a \leq k$,  $S^a_1,S^a_2,\dots,S^a_{p_a+q_a}\in \mbox{Fin}( \cup^\infty_{n=1} CU^{l_a,n})$ for some $m_a<l_a$.
\end{itemize}
\noindent
 such that 
\begin{itemize}
\item $\pi_{odd}(t_k)$ is in the orbit closure of $x$.
\item For all $1\leq a \leq k$ and $1\leq i \leq p_a+q_a$, $(\Theta^{l_a}_{t_k})^{-1}(S^a_i)\subset (\Theta^{m_a}_{t_k})^{-1}(S_a)$. 
\item For all $1\leq a \leq k$, the sets $\{(\Theta^{l_a}_{t_k})^{-1}(S^a_i)\}_{i=1}^{p_a+q_a}$ are disjoint.
\item For all $1\leq a \leq k$,  $(\Theta^{m_a}_{t_k})^{-1}(S_a)= \cup_{i=1}^{p_a} ((\Theta^{l_a}_{t_k})^{-1}(S^a_i))\gamma^a_i= \\ 
\cup_{j=p_a+1}^{p_a+q_a} ((\Theta^{l_a}_{t_k})^{-1}(S^a_j))\gamma^a_j \,.$
\end{itemize}
\noindent
Now, in the same way as in Step 0.  we construct a minimal landscape $t_{k+1}$ and elements $\gamma^{k+1}_1, \gamma^{k+1}_2,\dots, \gamma^{k+1}_{p_{k+1}+q_{k+1}}\in\Gamma$
and sets $S^{k+1}_1,S^{k+1}_2,\dots,S^{k+1}_{p_{k+1}+q_{k+1}}\in \mbox{Fin}( \cup^\infty_{n=1} CU^{l_{k+1},n})$ for some $m_{k+1}<l_{k+1}$ in such a way that
\begin{itemize}
\item $\pi_{odd}(t_{k+1})$ is in the orbit closure of $x$.
\item For all $1\leq a \leq k+1$ and $1\leq i \leq p_a+q_a$, $(\Theta^{l_a}_{t_{k+1}})^{-1}(S^a_i)\subset (\Theta^{m_a}_{t_{k+1}})^{-1}(S_a)$. 
\item For all $1\leq a \leq k+1$, the sets $\{(\Theta^{l_a}_{t_{k+1}})^{-1}(S^a_i)\}_{i=1}^{p_a+q_a}$ are disjoint.
\item For all $1\leq a \leq k+1$,  $(\Theta^{m_a}_{t_{k+1}})^{-1}(S_a)= \cup_{i=1}^{p_a} ((\Theta^{l_a}_{t_{k+1}})^{-1}(S^a_i))\gamma^a_i=\cup_{j=p_a+1}^{p_a+q_a} 
((\Theta^{l_a}_{t_{k+1}})^{-1}(S^a_j))\gamma^a_j \,.$
\end{itemize}
\noindent
Then we have a subsequence $k_1< k_2< \dots$ such that
$\lim_{r\to \infty} t_{k_r}=t\in A^\Gamma$ exists.
\begin{proposition} \label{utolso}
All $t$-local subset of $\Gamma$ is $t$-paradoxical.
\end{proposition}
\proof
Let $b\geq 1$ such that $(\Theta^{m_b}_t)^{-1}(S_b)$ is non-empty. 
We need to show that
\begin{enumerate}
\item For all  $1\leq i \leq p_b+q_b$, $(\Theta^{l_b}_{t})^{-1}(S^b_i)\subset (\Theta^{m_b}_{t})^{-1}(S_b)$. 
\item The sets $\{(\Theta^{l_b}_{t})^{-1}(S^b_i)\}_{i=1}^{p_b+q_b}$ are disjoint.
\item $(\Theta^{m_b}_{t})^{-1}(S_b)= \cup_{i=1}^{p_b} ((\Theta^{l_b}_{t})^{-1}(S^b_i))\gamma^b_i=\cup_{j=p_b+1}^{p_b+q_b} 
((\Theta^{l_b}_{t})^{-1}(S^b_j))\gamma^b_j \,.$
\end{enumerate}
\noindent
Now we  proceed in the same way as in the proof of Proposition \ref{clean}. Let us prove (1). 
Let $c>0$ be an integer such that $\gamma^b_1,\gamma^b_2\dots \gamma^b_{p_b+q_b}\in B_c(G,e_\gamma)$
Let $\gamma\in (\Theta^{l_b}_{t})^{-1}(S^b_i)$.  We need to show that
$\gamma\gamma^b_i\in  (\Theta^{m_b}_{t})^{-1}(S_b)$. Since $t=\lim_{r\to \infty} t_{k_r} $ we have $k_r>b$ such that
$\Theta^{c+l_b+m_b}_{t_{k_r}}(\gamma)=\Theta^{c+l_b+m_b}_t(\gamma)\,.$
Then $\gamma\in (\Theta^{l_b}_{t_{k_r}})^{-1}(S^b_i)$, so $\gamma\gamma^b_i\in (\Theta^{m_b}_{t_{k_r}})^{-1}(S_b)$, hence $\gamma\gamma^b_i\in (\Theta^{m_b}_{t})^{-1}(S_b)$.
The proof of (2) and (3) can be obtained in a similar fashion. \qed
\vskip 0.1in
\noindent
Observe that $\pi_{odd}(t)$ is in the orbit closure of $x$. Let $\hat{t}$ be a minimal landscape in the orbit closure of $t$ and let $Y$ be the totally finite part of the
orbit closure of  $\hat{t}$.
Then $\pi_{odd}(\hat{t})$ is in the orbit closure of $x$ as well. Hence, by Proposition \ref{freeminame}  the action of $\Gamma$ on $Y$ is free, minimal and amenable.
Also, by Proposition \ref{clean} all $\hat{t}$-local subsets of $\Gamma$ are $\hat{t}$-paradoxical. Hence, by Proposition \ref{key} the action of $\Gamma$ on $Y$ is
purely infinite. Since by Proposition \ref{legvege} $Y$ is homeomorphic to $\K$, our theorem follows. \qed

\section{A remark about actions on the compact Cantor set}
If $\Gamma$ is a non-amenable group then one can consider the compact Bernoulli
subshift $X=C^\Gamma$, where $C=\{0,1\}^\N$. Repeating the arguments of our 
paper one can construct a free, minimal purely infinite $\Gamma$-subshift $Y$ in 
$X$ such that $Y$ is homeomorphic to the Cantor set. If the group is exact, then using the witness-sets for Property A as in
Section \ref{river} one can even make the action amenable. This result is originally
due to R{\o}rdam and Sierakowski \cite{RS}. Our method just helps to avoid the
use of the \v{C}ech-Stone compactification. Similarly, one can eliminate the \v{C}ech-Stone compactification from the proof of Theorem 1.3. (iv) in \cite{Guent} and add
pure infinity to the properties of the action. That is, one can obtain  
(using the witness-sets of finite asymptotic dimension) the following result:
All countable non-amenable group $\Gamma$ of asymptotic dimension $d$ has
a free, minimal, purely infinite action of dynamic asymptotic
dimension at most $d$. One can also extend all these results for 
uniformly recurrent subgroups \cite{GW} as well.
Let $\Gamma$ be a finitely generated group and be $H\subset\Gamma$ a subgroup
such that the orbit closure of $H$ in $\Sub(\Gamma)$ is a closed, invariant,
minimal subspace. That is, $Z=\overline{O}(H)$ is a uniformly recurrent subgroup
(URS). If the Schreier graph $\Gamma/H$ is non-amenable (that is the URS Z is not
coamenable \cite{Elek}) then using the method of our paper we can construct
a free minimal purely infinite $Z$-proper (nonfree) action of $\Gamma$ 
(see \cite{Elek}
for the definition of $Z$-properness). If the Schreier graph is of
Property $A$ then we can even assume that the action is topologically amenable.

\end{document}